\documentclass[a4paper, 11pt]{amsart}
\usepackage[latin1]{inputenc}
\usepackage[T1]{fontenc}
\usepackage[english]{babel}
\usepackage{amssymb}
\usepackage{amsmath}
\usepackage{amsthm}
\usepackage{amscd}
\usepackage{amsfonts}
\usepackage{a4wide}
\usepackage{fancyhdr}
\usepackage{xypic}
\usepackage[shortlabels]{enumitem}

\newtheorem{ex}{Example}

\newcommand{\zn}[1]{\mathbb{Z}/#1\mathbb{Z}}
\newcommand{\F}{\mathbb{F}}
\newcommand{\Fp}{\mathbb{F}_{p^2}}
\newcommand{\etalchar}[1]{$^{#1}$}
\newcommand{\Mod}[1]{\ (\mathrm{mod}\ #1)}

\begin{document}

\title{On the supersingular GPST attack}
\author{{A}ndrea {B}asso and {F}abien {P}azuki}
\address{Andrea Basso. University of Birmingham, University Rd W, Birmingham B15 2TT, United Kingdom.}
\email{a.basso@cs.bham.ac.uk}
\address{Fabien Pazuki. University of Copenhagen, Institute of Mathematics, Universitetsparken 5, 2100 Copenhagen, Denmark, and Universit\'e de Bordeaux, IMB, 351, cours de la Lib\'eration, 33400 Talence, France.}
\email{fpazuki@math.ku.dk}
\thanks{Both authors thank Christophe Petit and Christophe Ritzenthaler for interesting feedback. The second author is a member of the project ANR-17-CE40-0012 Flair.}
\maketitle
\vspace{-0.9cm}
\begin{center}

\end{center}

\begin{abstract}
We explain why the first Galbraith-Petit-Shani-Ti attack on the Supersingular Isogeny Diffie-Hellman and the Supersingular Isogeny Key Encapsulation fails in some cases.
\end{abstract}

\vspace{0.5cm}

{\flushleft
\textbf{Keywords:} Isogenies, Supersingular elliptic curves, Modular invariants.\\
\textbf{Mathematics Subject Classification:} 14H52, 14K02, 11T71, 94A60, 81P94, 65P25. }

\vspace{0.5cm}

\begin{center}
---------
\end{center}

\vspace{1.5cm}

\thispagestyle{empty}

\section{Introduction}

In 2011, De Feo and Jao \cite{DFJ11} introduced Supersingular Isogeny Diffie-Hellman (SIDH), a post-quantum key exchange protocol that mimics the Diffie-Hellman protocol in the settings of isogenies between supersingular curves. In 2014, Jao, De Feo and Pl\^ut \cite{MR3259113} built upon SIDH to obtain a key encapsulation scheme called Supersingular Isogeny Key Encapsulation (SIKE).

We recently carried out a study of one of the most relevant attacks against SIDH and SIKE, namely the first attack presented in \cite{Galbraith2016OnCryptosystems}. We call it the Galbraith-Petit-Shani-Ti attack, or GPST. It is an active attack where the attacker impersonates one of the two parties and can recover the static key of the other party in about as many interactions as the number of bits of the key. We will show in Section \ref{failure} that the attack may fail in some precise circumstances.

\subsection*{Preliminaries}

For an elliptic curve $E$ defined over a field $k$, we denote the base point of the elliptic curve by $\mathcal{O}$ and the modular invariant of $E$ by $j(E)$. The invariant is an element of $k$ that characterizes the $\overline{k}$-isomorphism class of $E$. For $P\in{E(k)}$ rational over $k$, we denote by $\langle P \rangle$ the subgroup of $E(k)$ generated by $P$.

The parameters of the SIDH protocol, in the form and notation of \cite{Galbraith2016OnCryptosystems}, are the following:
\begin{itemize}
    \item A prime $p = 2^n3^mf - 1$, where $n, m, f$ are natural numbers, $f$ is small and $2^n \approx 3^m$.
    \item A supersingular elliptic curve $E_0$ defined over $\Fp$.
    \item Points $P_A, Q_A \in E_0$ which form a basis of $E_0[2^n]$ and points $P_B, Q_B \in E_0$ which form a basis of $E_0[3^m]$.
\end{itemize}

We refer the reader to the original paper \cite{Galbraith2016OnCryptosystems} for the presentation of the attack. We introduce here the two assumptions that are required for the GPST attack to take place:
\begin{enumerate}
    \item[(a)] Alice uses a static key $(a_1, a_2)$. The values $a_1, a_2$ are elements of $\zn{2^n}$, not both divisible by 2.
    \item[(b)] The attacker has access to an oracle $O$ such that, if $E, E'$ denote supersingular elliptic curves defined over $\Fp$ and $R, S$ denote any two points on $E$,
        \begin{equation}\label{oracle}
            O(E, R, S, E') =
                \begin{cases}
                    true, &\text{if } j(E / \langle[a_1]R + [a_2]S\rangle) = j(E'),\\
                    false, &\text{otherwise.}
                \end{cases}
        \end{equation}
        Such an oracle can be realized in practice by having a mechanism in place that ensures the two parties obtain the same shared key.
\end{enumerate}
The attacker does not need any additional information, which makes the GPST attack one of the most powerful attacks known in the literature.

\section{Attacking the attack}\label{failure}
Alice's static key $(a_1, a_2)$ is always equivalent (\textit{i.e.} it leads to the same key exchange) to a key of the form $(1, \alpha)$ or $(\alpha, 1)$, where $\alpha$ is again an element of $\zn{2^n}$ \cite[Lemma 2.1]{Galbraith2016OnCryptosystems}. Without loss of generality, we may assume we are in the former case.

Each iteration of the GPST attack relies on the following implication, used in the paragraph \textit{First step of the attack} of \cite{Galbraith2016OnCryptosystems}:
\begin{equation}\label{wrong}
\Big(j(E/\langle T_1 \rangle) = j(E/\langle T_2 \rangle) \Big) \Rightarrow \Big( \langle T_1 \rangle = \langle T_2 \rangle\Big)
\end{equation}
to recover a single bit of the key. In particular, the oracle query (\ref{oracle}) is used to recover whether $\langle R' + \alpha S' \rangle$ is equal to $\langle R + \alpha S\rangle$, where $R, S, R', S'$ are points on an elliptic curve computed by the attacker. The points $R'$ and $S'$ are chosen such that $\langle R' + \alpha S' \rangle = \langle R + \alpha S\rangle$ implies the $i$-th bit $\alpha_i$ of the key is a zero-bit. Thus, if the oracle returns \textit{true}, then $\alpha_i=0$, otherwise $\alpha_i=1$.

However, implication (\ref{wrong}) is not correct, because there exist pairs of elliptic curves $E_1, E_2$ with multiple isogenies between them. Thus, if $\phi_1$ and $\phi_2$ are distinct cyclic isogenies between $E_1$ and $E_2$ and $K_1$ and $K_2$ are their respective kernels, we have that $E_1/K_1 $ is isomorphic to $E_1/K_2 $ because they are both isomorphic to $E_2$, but $K_1$ and $K_2$ may not be equal and need not even be isomorphic. Over $\mathbb{Q}$, a classical example can be found in \cite[page 110]{Sil94}, namely the pair $(E, E)$, where $E$ is the curve defined over $\mathbb{Q}$ with $j(E)=8000$, which is given by the affine Weierstrass equation $y^2=x^3+4x^2+2x$. This curve $E$ admits a rational $2$-isogeny $\varphi$ to itself (hence cyclic, and given explicitly in \textit{loc. cit.}), and one computes $j(E/\ker\varphi)=8000=j(E)$, hence $E$ and $E/\ker\varphi$ are in fact $\overline{\mathbb{Q}}$-isomorphic, but clearly $\ker\varphi\neq\{O\}$.

Over a finite field of characteristic $p>0$, if one starts with a supersingular elliptic curve $E$, then for any prime $\ell\neq p$, for any positive integer $m$ and for any cyclic subgroup $G$ of order $\ell^m$, the elliptic curve $E/G$ will also be supersingular, because there are still no points of order $p$ on $E/G$. As there are only finitely many supersingular curves over $\overline{\mathbb{F}}_p$ (see for instance \cite[Theorem 4.1, pages 148-149]{Sil86}), there exist pairs of cyclic groups $(G_1, G_2)$ where $G_1$ and $G_2$ are not isomorphic and such that $E/G_1$ is isomorphic to $E/G_2$. This argument is used to compute the endomorphism ring of supersingular elliptic curves in \cite[page 146]{Sil86} (see also \cite[page 267]{Hus04}).

Let us give now a detailed example where the Galbraith-Petit-Shani-Ti attack fails to recover the private key of a SIDH key exchange.

\begin{ex}\label{wrongexample}
Let $p = 2^5 3^3 - 1 = 863$ be a prime and $k$ be the finite field $\mathbb{F}_{p^2}$, considered as $\F_p(\beta)$, where $\beta$ satisfies the quadratic equation $\beta^2 -\beta + 5 = 0$. Let $E_0$ be the supersingular elliptic curve with affine Weierstrass model $y^2 = x^3 + (531\beta + 538)x + (720\beta +375)$ over $k$. Consider the points 
\begin{gather*}
    P_A = (834\beta + 726, 642\beta + 130), \quad Q_A = (583\beta + 276, 180\beta + 854), \\
    P_B = (254\beta + 697, 516\beta + 268), \quad Q_B = (753\beta + 317, 234\beta + 532).
\end{gather*}
The points $P_A, Q_A$ form a basis of $E_0[2^5]$, and $P_B, Q_B$ form a basis of $E_0[3^3]$.
Let Alice's key be of the form $(1, \alpha)$, with $\alpha = 10$ (written in base ten). The value $\alpha$ can also be expressed in bits as $1010$ (written in base two). From now on, by a slight abuse of notation we refer to $\alpha$ as the key.

Let us now carry on the Galbraith-Petit-Shani-Ti attack. Assume the randomly generated values $b_1, b_2$ are $b_1 = 1, b_2 = 6$ (for simplicity, we assume $b_1$ and $b_2$ stay constant across iterations, but we only need $b_1 = 1, b_2 = 6$ in the second round for the attack to fail).
Let $\phi_A$ be the isogeny with kernel $K_A = \langle P_A + [\alpha]Q_A\rangle$. Then $E_A = E_0/K_A$ has affine Weierstrass model $y^2 = x^3 + (40\beta+535)x + (720\beta+768)$.
Let us begin the attack and let $\phi_B$ be the isogeny with kernel $K_B = \langle P_B + [6]Q_B\rangle$. Thus $E_B = E_0/K_B$ has affine Weierstrass model $y^2 = x^3 + 105x + 254$. Furthermore, we have
\begin{gather*}
    R = \phi_B(P_A) = (151\beta + 257, 594\beta + 2), \\
    S = \phi_B(Q_A) = (98\beta + 386, 286\beta + 58).
\end{gather*}
Note that since the degree of $\phi_B$ is coprime with the order of $P_A$ and $Q_A$, the points $R$ and $S$ have also order $2^5$.

The first iteration of the attack starts by computing
\begin{gather*}
    \theta_0 = \sqrt{(1 + 2^4)^{-1}}  \Mod{2^5} = 7, \\
    R'_0 = [\theta_0]R = (527\beta + 129, 700\beta + 163), \\
    S'_0 = [\theta_0][1 + 2^4]S = (164\beta + 377, 566\beta + 641).
\end{gather*}
It then proceeds by querying the oracle with $O(E_B, R'_0, S'_0, E_{AB}) = E_B/\langle \phi_A(P_B) + [6]\phi_A(Q_B) \rangle$. The curve $E_{AB}: y^2 = x^3 + 698x + (516\beta+605)$ has $j$-invariant $117$. The curve $E_B / \langle R'_0 + [\alpha]S'_0 \rangle$ also has $j$-invariant $117$, thus the oracle response is \emph{true} and the first (rightmost) bit of the key is a zero-bit. Hence the attack correctly obtains the first bit of the key.

For the second round, the attacker computes
\begin{gather*}
    \theta_1 = \sqrt{(1 + 2^3)^{-1}} \Mod{2^5} = 5, \\
    R'_1 = [\theta_1]R = (97\beta+ 261, 795\beta + 545), \\
    S'_1 = [\theta_1][1 + 2^3]S = (718\beta + 214, 450\beta + 844),
\end{gather*}
and queries the oracle on $O(E_B, R'_1, S'_1, E_{AB})$. As before, the curves $E_{AB}$ and $E_B / \langle R'_1 + [\alpha]S'_1 \rangle$ both have $j$-invariant $117$, thus the oracle responds \emph{true}. Hence this time the attacker incorrectly deduces the second bit (reading from right to left) of the key to be zero.

Brute-forcing the remaining bits does not yield any solution, since there is no key $\alpha$ ending with two zero-bits (00) such that $E_0/\langle P_A + [\alpha] Q_A \rangle$ is isomorphic to $E_A$.\\

\end{ex}

The previous example establishes that there exist cases where the attack fails. The following are sufficient conditions (when considered together) for the attack to fail at the $i$-th iteration of the attack, for $1 \leq i <  n - 3$ (the attack brute-forces the last three bits of the key and thus cannot fail after $n-3$):

\begin{enumerate}[(i)]
    \item There exist two distinct isogenies $\phi_1$ and $\phi_2$ between $E_B$ and $E_{AB}$. \label{firstcond}
    \item One has $R + [\alpha]S = P_1$ and $R'_i + [\alpha]S'_i = P_2$, where $P_1$ and $P_2$ are generators of the kernel of $\phi_1$ and $\phi_2$, respectively. \label{secondcond}

    \item The $i$-th key bit $\alpha_i$ is $1$.
    \label{thirdcond}
\end{enumerate}

Determining whether condition \ref{firstcond} is satisfied is computationally hard since it is generally hard to compute isogenies between two given curves. Moreover, the points $P_1$ and $P_2$ are dependent on the key $\alpha$ and thus unknown to the attacker. Hence, while the attacker influences the points $P_1$ and $P_2$ by choosing $\phi_B$, they cannot know whether his choice of $\phi_B$ would cause the attack to fail.

Condition \ref{secondcond} ensures that the attack is taking place in the case where the points $R$ and $S$ give rise to the isogenies considered in the first condition. Given the points $P_1, P_2$, we have
\begin{equation*}
	\begin{cases}
		P_1 = R + [\alpha]S,\\
		P_2 = R'_i + [\alpha]S'_i = [\theta_i]R + [\theta_i][\alpha + \alpha_i2^{n - i - 1}]S,
	\end{cases}
\end{equation*}
thus $[\theta_i]P_1 - P_2 = [\theta_i][2^{n-i-1}][\alpha_i]S = [\theta_i][2^{n-i-1}]S$ (because of condition \ref{thirdcond}). If a $\theta_i$ exists such that the previous equation is satisfied, then suitable points $R$ and $S$ also exist.

Condition \ref{thirdcond} is necessary because the attack fails by incorrectly deducing a zero-bit instead of a one-bit.

These conditions appear to be quite rare. If it could be shown that the number of cases where the attack fails is polylogarithmic in the security parameter, the attacker could simply retry the attack with different choices of $\phi_B$ until it succeeds. The attack would still be efficient.\\

Let us add a remark. Each iteration of the attack depends on the previous one since the modified value $R'_i$ of iteration $i+1$ relies on $K_i$, whose last bit is obtained in the $i$\textsuperscript{th} iteration. To see how the error propagates, assume the attack fails at iteration $i$. Thus $\alpha_i$, the $i$\textsuperscript{th} bit of the key, is a one-bit but the attacker deduces it to be a zero-bit. It follows that $K_i$ = $K_{i-1}$ for the attacker since prepending zero-bits does not affect the value of $K$. The attacker then computes the points
\begin{equation*}
    R'_i = [\theta]R - [\theta][2^{n-i-1}K_{i-1}]S, \quad S'_i = [\theta][1+ 2^{n-i-1}]S
\end{equation*}
and queries the oracle on $O(E_B, R'_i, S'_i, E_{AB})$.

Thus, omitting $\theta$ since it does not affect the subgroup, we have
\begin{align*}
    &\langle R'_i + [\alpha]S'_i \rangle \\
    = ~ &\langle R + [-2^{n -i -1}K_{i-1} + \alpha (1 + 2^{n-i-1})]S \rangle \\
    = ~ &\langle R + [\alpha]S + [-2^{n -i -1}K_{i-1} + 2^{n-i-1}(K_{i-1} + 2^{i-1} + 2^i\alpha_i + 2^{i+1}\alpha')]S \rangle \\
    = ~ &\langle R + [\alpha]S + [2^{n-2} + 2^{n-1}\alpha_i]S \rangle \\
    \neq ~ &
      \langle R + [\alpha]S \rangle,
\end{align*}
To see why the the inequality in the last line holds, consider that if the two cyclic subgroups are the same, their generators must be multiples of each other. Assume then $R + [\alpha]S + [2^{n-2} + 2^{n-1}\alpha_i]S = [m](R + [\alpha]S)$. The linear independence of $R$ and $S$ implies that $R = [m]R$, which means $m \equiv 1 \bmod 2^n$. Thus, the equality on $S$ implies that $ [2^{n-2} + 2^{n-1}\alpha_i]S = \mathcal{O}$, which is impossible because $\alpha_i \in \{0, 1\}$ and $0 \not\equiv 2^{n-2} \not\equiv -2^{n-1} \bmod 2^n$ . This means that the oracle will return \textit{false} (unless the two subgroups, $\langle R + [\alpha]S \rangle$ and $\langle R'_i + [\alpha]S'_i \rangle$, give rise again to two isogenies with isomorphic codomain).

The subsequent iterations further propagate the error and the oracle will respond \textit{false} to every further query (unless, again, the isogenies with different kernels have the same codomain, which is believed to be a rare occurrence). Thus, if the attack fails at the $i$\textsuperscript{th} iteration, all the following key bits are deduced to be a one-bit. This allows the attacker to approximately identify the part of the key that has been correctly deduced and target, with a different choice of $\phi_B$, only the remaining part.

\section{Sage Implementation}

In this section, we report the source code of a SAGE \cite{sage} implementation that shows the attack failing in the case detailed in Example \ref{wrongexample}.

\begin{verbatim}

def oracle(E, R, S, Ep):
    EAp = E.isogeny(R + alpha*S).codomain()
    return EAp.j_invariant() == Ep.j_invariant()

## Performs the GPST attack, assuming that Alice's key
## is of the form (1, alpha).
def attack(n, m, E0, PA, QA, PB, QB, EA, phiAPB, phiAQB):
    K = 0

    for i in range(n - 3):
        alpha = 0
        
        b1 = 1
        b2 = 6
        
        KB = b1*PB + b2*QB
        phiB = E0.isogeny(KB)
        EB = phiB.codomain()
        
        R = phiB(PA)
        S = phiB(QA)
        EAB = EA.isogeny(b1*phiAPB + b2*phiAQB).codomain()

        FF = IntegerModRing(2^n)
        theta = Integer(FF((1 + 2^(n - i - 1))^-1).sqrt())
        Rprime = theta * (R - (2^(n - i - 1) * K) * S)
        Sprime = theta * (1 + 2^(n - i - 1)) * S

        response = oracle(EB, Rprime, Sprime, EAB)

        if response == False:
            alpha = 1

        K += alpha*2^i

    found = False

    ## Bruteforcing the rest of the key
    for i in range(2):
        for j in range(2):
            for k in range(2):
                key = K
                key += i*2^(n-3) + j*2^(n-2) + k*2^(n-1)

                EAprime = E0.isogeny(PA + key*QA).codomain()

                if EAprime.j_invariant() == EA.j_invariant():
                    solution = (1, key % 2^n)
                    found = True
                    break

    if found:
        return solution
    else:
        return "Key not found"

## Setup
lA = 2
lB = 3
eA = 5
eB = 3    
f = 1

p = lA^eA*lB^eB*f - 1 #p in Primes() returns True
F.<x> = GF(p^2)

# E0: y^2 = x^3 + (531\beta + 538)x + 720\beta + 375
# E0 is supersingular, as shown by E0.is_supersingular()
E0 = EllipticCurve(F, [531*x + 538, 720*x + 375])

PA = E0(834*x + 726, 642*x + 130)
QA = E0(583*x + 276, 180*x + 854) #PA and QA form a basis of E0[2^5]
PB = E0(254*x + 697, 516*x + 268)
QB = E0(753*x + 317, 234*x + 532) #PB and QB form a basis of E0[3^3]

mA = 1
nA = alpha = 10

KA = mA*PA + nA*QA

phiA = E0.isogeny(KA)
EA = phiA.codomain()

phiAPB = phiA(PB)
phiAQB = phiA(QB)

alphaprime = attack(eA, eB, E0, PA, QA, PB, QB, EA, phiAPB, phiAQB)
print("The key is: %s" % str(alphaprime))
\end{verbatim}

\end{document}